\documentclass[11pt]{article}
\usepackage{mathrsfs}
\usepackage{amsthm}
\usepackage{amssymb}
\usepackage{amsmath}
\usepackage{graphicx}
\usepackage{color}
\usepackage{amsfonts}
\usepackage{float}
\usepackage{cite}
\usepackage[text={140mm,210mm},left=35mm,vmarginratio=1:1]{geometry}
\newtheorem{theorem}{Theorem}[section]

\newtheorem{lemma}[theorem]{Lemma}
\newtheorem{proposition}[theorem]{Proposition}

\numberwithin{equation}{section}
\normalsize

\begin{document}
\title{\textbf{The critical infection rate of the high-dimensional two-stage contact process}}

\author{Xiaofeng Xue \thanks{\textbf{E-mail}: xfxue@bjtu.edu.cn \textbf{Address}: School of Science, Beijing Jiaotong University, Beijing 100044, China.}\\ Beijing Jiaotong University}

\date{}
\maketitle

\noindent {\bf Abstract:} In this paper we are concerned with the two-stage contact process on the lattice $\mathbb{Z}^d$ introduced in \cite{Krone1999}. We gives a limit theorem of the critical infection rate of the process as the dimension $d$ of the lattice grows to infinity. A linear system and a two-stage SIR model are two main tools for the proof of our main result.

\quad


\section{Introduction}\label{section one}

In this paper we are concerned with the two-stage contact process on the lattice $\mathbb{Z}^d$, which is introduced in \cite{Krone1999}. For each $x\in \mathbb{Z}^d$, we use $\|x\|$ to denote the $l_1$-norm of $x$, i.e.,
\[
\|x\|=\sum_{i=1}^d|x_i|
\]
for $x=(x_1,\ldots,x_d)$. For any $x,y\in \mathbb{Z}^d$, we write $x\sim y$ when and only when $\|x-y\|=1$. In other words, we use $x\sim y$ to denote that $x$ and $y$ are neighbors. We use $O$ to denote the origin $(0,0,\ldots,0)$.

The two-stage contact process $\{\eta_t\}_{t\geq 0}$ on $\mathbb{Z}^d$ is a continuous time Markov process with state space $X_1=\{0,1,2\}^{\mathbb{Z}^d}$. The transition rates function is given as follows.
\begin{equation}\label{equ 1.1 transition rate}
\eta_t(x)\rightarrow i \text{~at rate~}
\begin{cases}
&1 \text{~if~} \eta_t(x)=2 \text{~and~} i=0,\\
&\gamma \text{~if~} \eta_t(x)=1 \text{~and~} i=2,\\
&1+\delta \text{~if~} \eta_t(x)=1 \text{~and~} i=0,\\
&\lambda\sum_{y:y\sim x}1_{\{\eta_t(y)=2\}} \text{~if~} \eta_t(x)=0 \text{~and~} i=1,\\
&0 \text{~otherwise}
\end{cases}
\end{equation}
for each $x\in Z^d$ and $t\geq 0$, where $\lambda,\gamma,\delta$ are positive constants and $1_A$ is the indicator function of the event $A$. The constant $\lambda$ is called the infection rate of the process.

The process $\{\eta_t\}_{t\geq 0}$ intuitively describes the spread of an epidemic on $\mathbb{Z}^d$. Each vertex is in one of three states, which are `healthy', `semi-infected'
and `fully-infected'. A healthy vertex is infected at rate proportional to the number of fully-infected neighbors to become a `semi-infected' one while a `semi-infected' vertex waits for an exponential time with rate $\gamma$ to become fully-infected or waits for an exponential time with rate $1+\delta$ to become healthy. A `fully-infected' vertex waits for an exponential time with rate $1$ to become healthy.

The two-stage contact process $\{\eta_t\}_{t\geq 0}$ is introduced by Krone in \cite{Krone1999}, where a duality relation between this two-stage contact process and an `on-off' process is given. In \cite{Fox2015}, Fox gives a simple proof of the duality relation given in \cite{Krone1999} and answers most of the open questions posed in \cite{Krone1999}.

If $\gamma=+\infty$, i.e., a semi-infected vertex becomes a fully-infected one immediately, then equivalently there is only one infected state for the process and hence the model reduces to the classic contact process introduced in \cite{Har1974}. For a survey of the classic contact process, see Chapter 6 of \cite{Lig1985} and Part \uppercase\expandafter{\romannumeral1} of \cite{Lig1999}.

\section{Main result}\label{section two}
In this section we give our main result. First we introduce some notations and definitions. Throughout this paper we assume that $\{x:\eta_0(x)=1\}=\emptyset$, i.e., there is no semi-infected vertex at $t=0$.  For each $A\subseteq \mathbb{Z}^d$, we write $\eta_t$ as $\eta_t^A$
when $\{x:\eta_0(x)=2\}=A$. If $A=\{x\}$ for some $x\in \mathbb{Z}^d$, then we write $\eta_t^A$ as $\eta_t^x$ instead of $\eta_t^{\{x\}}$. When we omit the superscript $A$, then we mean that $A=\mathbb{Z}^d$. For any $t\geq 0$, we use $C_t^A$ to denote
\[
\big\{x:\eta_t^A(x)=2\big\}
\]
as the set of fully-infected vertices at the moment $t$. We denote by $P_d^{\lambda,\gamma,\delta}$ the probability measure of the two-stage contact process $\{\eta_t\}_{t\geq 0}$ with parameter $\lambda,\gamma,\delta$ defined as in Equation \eqref{equ 1.1 transition rate}. It is obviously that $P_d^{\lambda,\gamma,\delta}\big(C_t^O\neq \emptyset \text{~for all~}t\geq 0\big)$ is increasing with $\lambda$. Hence it is reasonable to define
\begin{equation}\label{equ 2.1}
\lambda_c(d,\gamma,\delta)=\sup\Big\{\lambda:P_d^{\lambda,\gamma,\delta}\big(C_t^O\neq \emptyset \text{~for all~}t\geq 0\big)=0\Big\}.
\end{equation}
$\lambda_c(d,\gamma,\delta)$ is called the critical infection rate of the two-stage contact process, with infection rate below which fully-infected vertices die out with probability one conditioned on $O$ is the unique fully-infected vertex at $t=0$. Now we give our main result, which is a limit theorem of $\lambda_c(d,\gamma,\delta)$ as the dimension $d$ grows to infinity.
\begin{theorem}\label{theorem main 2.1}
For any $\gamma,\delta>0$, if $\lambda_c(d,\gamma,\delta)$ is defined as in Equation \eqref{equ 2.1}, then
\[
\lim_{d\rightarrow+\infty}2d\lambda_c(d,\gamma,\delta)=1+\frac{1+\delta}{\gamma}.
\]
\end{theorem}

\proof[Remark 1]
Let $\alpha_c(d)$ be the critical infection rate for the classic contact process on $\mathbb{Z}^d$, then it is shown in \cite{Grif1983} that
\begin{equation}\label{equ 2.2}
\lim_{d\rightarrow+\infty}2d\alpha_c(d)=1.
\end{equation}
Our main result can be considered as an extension of conclusion \eqref{equ 2.2} since when $\gamma=+\infty$ the two-stage contact process reduces to the classic contact process.

\qed

\proof[Remark 2]

It is shown in \cite{Fox2015} that $\lambda_c(d,\gamma,\delta)=+\infty$ when $\gamma<\frac{1}{4d-1}$. Our main result is not in contrast with this conclusion since $\gamma>\frac{1}{4d-1}$ for sufficiently large $d$.

\qed

\proof[Remark 3]

It is shown in \cite{Grif1983} that $\alpha_c(d)\leq \frac{1}{2d}+\frac{1}{2d^2}+o(\frac{1}{d^2})$. Hence it is natural to guess that there exists $f(\gamma,\delta)>0$ such that
\[
\lambda_c(d,\gamma,\delta)\leq \frac{1}{2d}(1+\frac{1+\delta}{\gamma})+\frac{f(\gamma,\delta)}{d^2}+o(\frac{1}{d^2}).
\]
However, according to our current approach we have not managed to obtain such a $f$ yet. We will work on this question as a further study.

\qed

The proof of Theorem \ref{theorem main 2.1} is divided into Section \ref{section three} and Section \ref{section four}. In Section \ref{section three}, we give the proof of $\liminf_{d\rightarrow+\infty}2d\lambda_c(d,\gamma,\delta)\geq 1+\frac{1+\delta}{\gamma}$. For this purpose, we will introduce a linear system with state space $\{\mathbb{Z}^2_+\}^{\mathbb{Z}^d}$ as a main auxiliary model, where $\mathbb{Z}_+=\{0,1,2,\ldots\}$. In Section \ref{section four}, we give the proof of $\limsup_{d\rightarrow+\infty}2d\lambda_c(d,\gamma,\delta)\leq 1+\frac{1+\delta}{\gamma}$. The proof is inspired by the approach introduced in \cite{Xue2017}. We will introduce a two-stage SIR(susceptible-infected-recovered) model, the critical infection rate of which is an upper bound of $\lambda_c(d,\gamma,\delta)$.

\section{The proof of $\liminf_{d\rightarrow+\infty}2d\lambda_c(d,\gamma,\delta)\geq 1+\frac{1+\delta}{\gamma}$}\label{section three}
In this section we give the proof of $\liminf_{d\rightarrow+\infty}2d\lambda_c(d,\gamma,\delta)\geq 1+\frac{1+\delta}{\gamma}$. First we introduce an auxiliary model, which is a linear system with state space $\{\mathbb{Z}_+^2\}^{\mathbb{Z}^d}$. For a survey of the linear system, see Chapter 9 of \cite{Lig1985}. Let
$\{(\zeta_t,\theta_t)\}_{t\geq 0}$ be a continuous-time Markov process with state space $\{\mathbb{Z}_+^2\}^{\mathbb{Z}^d}$, where $\mathbb{Z}_+=\{0,1,2,\ldots\}$. That is to say, at each vertex $x$ there is an vector $\big(\zeta(x),\theta(x)\big)$. The transition rates function of $\{(\zeta_t,\theta_t)\}_{t\geq 0}$ is given as follows. For each $x\in \mathbb{Z}^d$ and $t\geq 0$,
\begin{align}\label{equ 3.1 transition rate of linear system}
&\big(\zeta_t(x),\theta_t(x)\big)\rightarrow (a,b) \text{~at rate~}\\
&\begin{cases}
1 \text{~if~} a=b=0,\\
\delta \text{~if~} a=\zeta_t(x)\text{~and~}b=0,\\
\gamma \text{~if~} a=\zeta_t(x)+\theta_t(x)\text{~and~}b=0,\\
\lambda \text{~if~} y\sim x, a=\zeta_t(x)\text{~and~}b=\theta_t(x)+\zeta_t(y),\\
0 \text{~otherwise}.
\end{cases}
\notag
\end{align}
The auxiliary model $\{(\zeta_t,\theta_t)\}_{t\geq 0}$ and the two-stage contact process have the following coupling relationship.
\begin{lemma}\label{lemma 3.1}
For $x\in \mathbb{Z}^d$ and $t\geq 0$, let
\[
\widehat{\eta}_t(x)=
\begin{cases}
2 \text{~if~} \zeta_t(x)>0,\\
1 \text{~if~} \zeta_t(x)=0 \text{~and~}\theta_t(x)>0,\\
0 \text{~if~} \zeta_t(x)=\theta_t(x)=0,
\end{cases}
\]
then $\{\widehat{\eta}_t\}_{t\geq 0}$ is a two-stage contact process with transition rates function given in Equation \eqref{equ 1.1 transition rate}.
\end{lemma}

\proof[Proof of Lemma \ref{lemma 3.1}]
According to Equation \eqref{equ 3.1 transition rate of linear system} and the definition of $\widehat{\eta}_t$,
$\widehat{\eta}_t(x)$ flips from $2$ to $0$ when and only when $\zeta_t(x)>0$ and $(\zeta_t(x),\theta_t(x))$ flips to $(0,0)$, the transition rate of which is $1$. $\widehat{\eta}_t(x)$ flips from $1$ to $0$ when and only when $\zeta_t(x)=0$ and $\theta_t(x)$ flips from a positive state to $0$, the transition rate of which is $1+\delta$. $\widehat{\eta}_t(x)$ flips from $1$ to $2$ when and only when $\zeta_t(x)=0,\theta_t(x)>0$ and $\big(\zeta_t(x),\theta_t(x)\big)$ flips to
\[
\big(\zeta_t(x)+\theta_t(x),0\big)=\big(\theta_t(x),0\big),
\]
the transition rate of which is $\gamma$.
$\widehat{\eta}_t(x)$ flips from $0$ to $1$ when and only when $\zeta_t(x)=\theta_t(x)=0$ and $\big(\zeta_t(x),\theta_t(x)\big)$ flips from $(0,0)$ to
\[
\big(\zeta_t(x),\theta_t(x)+\zeta_t(y)\big)=\big(0,\zeta_t(y)\big)
\]
for some $y\sim x$ that $\zeta_t(y)>0$, the transition rate of which is
\[
\lambda\sum_{y:y\sim x}1_{\{\zeta_t(y)>0\}}=\lambda\sum_{y:y\sim x}1_{\{\widehat{\eta}_t(y)=2\}}.
\]
In conclusion, the transition rates function of $\{\widehat{\eta}_t\}_{t\geq 0}$  is as that given in Equation \eqref{equ 1.1 transition rate}.

\qed

According to Lemma \ref{lemma 3.1}, we can consider that the two-stage contact process $\{\eta_t\}_{t\geq 0}$ and $\{(\zeta_t,\theta_t)\}_{t\geq 0}$ are defined under the same probability space that
\[
\eta_t(x)=2\times 1_{\{\zeta_t(x)>0\}}+1_{\{\zeta_t(x)=0\text{~and~}\theta_t(x)>0\}}
\]
for each $x\in \mathbb{Z}^d$.

Conditioned on all the vertices are in state $2$ at $t=0$, it is shown in \cite{Krone1999} that the distribution of $\eta_t$ converges weakly to a probability distribution $\nu=\nu^{\lambda,\gamma,\delta}$ on $\{0,1,2\}^{\mathbb{Z}^d}$ as $t$ grows to infinity. As a result,
\begin{equation}\label{equ 3.2}
\lim_{t\rightarrow+\infty}P^{\lambda,\gamma,\delta}_d\big(\eta_t(O)=1\text{~or~}2\big)=\nu^{\lambda,\gamma,\delta}\big(\eta(O)=1\text{~or~}2\big).
\end{equation}
Note that when we omit the superscript of $\eta_t$ we mean that all the vertices are in state $2$ at $t=0$.

Our proof of $\liminf_{d\rightarrow+\infty}2d\lambda_c(d,\gamma,\delta)\geq 1+\frac{1+\delta}{\gamma}$ relies heavily on the following proposition, which is Theorem 1.2 of \cite{Fox2015}.

\begin{proposition}\label{proposition 3.2}
(Fox, 2015) $P^{\lambda,\gamma,\delta}\big(C_t^O\neq \emptyset\text{~for all~}t\geq 0\big)>0$ if and only if
\[
\nu^{\lambda,\gamma,\delta}\big(\eta(O)=1\text{~or~}2\big)>0.
\]
\end{proposition}
For the proof of this proposition, see Section 3.5 of \cite{Fox2015}. Now we give the proof of $\liminf_{d\rightarrow+\infty}2d\lambda_c(d,\gamma,\delta)\geq 1+\frac{1+\delta}{\gamma}$.

\proof [Proof of $\liminf_{d\rightarrow+\infty}2d\lambda_c(d,\gamma,\delta)\geq 1+\frac{1+\delta}{\gamma}$]

We assume that $\zeta_0(x)=1$ and $\theta_0(x)=0$ for all $x\in \mathbb{Z}^d$, then according to Lemma \ref{lemma 3.1} and Markov's inequality,
\begin{align}\label{equ 3.3}
P_d^{\lambda,\gamma,\delta}\big(\eta_t(O)=1\text{~or~}2\big)&=P_d^{\lambda,\gamma,\delta}\big(\zeta_t(O)\geq 1\big)+P_d^{\lambda,\gamma,\delta}\big(\zeta_t(O)=0,\theta_t(O)\geq 1\big)\notag\\
&\leq P_d^{\lambda,\gamma,\delta}\big(\zeta_t(O)\geq 1\big)+P_d^{\lambda,\gamma,\delta}\big(\theta_t(O)\geq 1\big) \notag\\
&\leq E_d^{\lambda,\gamma,\delta}\zeta_t(O)+E_d^{\lambda,\gamma,\delta}\theta_t(O),
\end{align}
where $E_d^{\lambda,\gamma,\delta}$ is the expectation operator with respect to $P_d^{\lambda,\gamma,\delta}$. According to the transition rates function of $(\zeta_t,\theta_t)$ and Theorem 9.1.27 of \cite{Lig1985}, which is an extended version of Hille-Yosida Theorem for the linear system,
\begin{equation*}
\begin{cases}
&\frac{d}{dt}E_d^{\lambda,\gamma,\delta}\zeta_t(O)=-E_d^{\lambda,\gamma,\delta}\zeta_t(O)+\gamma E_d^{\lambda,\gamma,\delta}\theta_t(O),\\
&\frac{d}{dt}E_d^{\lambda,\gamma,\delta}\theta_t(O)=-(1+\gamma+\delta)E_d^{\lambda,\gamma,\delta}\theta_t(O)+\lambda\sum_{y:y\sim O}E_d^{\lambda,\gamma,\delta}\zeta_t(y).
\end{cases}
\end{equation*}
Conditioned on $\zeta_0(x)=1, \theta_0(x)=0$ for all $x\in \mathbb{Z}^d$, $E_d^{\lambda,\gamma,\delta}\zeta_t(y)$ does not depend on the choice of $y$ according to the spatial homogeneity of our process, hence
\begin{equation}\label{equ 3.4}
\frac{d}{dt}
\begin{pmatrix}
E_d^{\lambda,\gamma,\delta}\zeta_t(O)\\
E_d^{\lambda,\gamma,\delta}\theta_t(O)
\end{pmatrix}
=
\begin{pmatrix}
-1&&\gamma\\
&&\\
2d\lambda && -(1+\gamma+\delta)
\end{pmatrix}
\begin{pmatrix}
E_d^{\lambda,\gamma,\delta}\zeta_t(O)\\
E_d^{\lambda,\gamma,\delta}\theta_t(O)
\end{pmatrix}.
\end{equation}
We use $G$ to denote
\[
\begin{pmatrix}
-1&&\gamma\\
&&\\
2d\lambda && -(1+\gamma+\delta)
\end{pmatrix}.
\]
Let $c_1,c_2$ be the two eigenvalues of $G$, then $E_d^{\lambda,\gamma,\delta}\zeta_t(O)=a_1e^{c_1t}+a_2e^{c_2t}$ and
$E_d^{\lambda,\gamma,\delta}\theta_t(O)=b_1e^{c_1t}+b_2e^{c_2t}$ for any $t\geq 0$ according to Equation \eqref{equ 3.4}, where $a_1,a_2,b_1,b_2$ are four constants. When $2d\lambda\gamma<1+\gamma+\delta$, it is easy to check that ${\rm Re}(c_1), {\rm Re}(c_2)<0$ and hence
\[
\lim_{t\rightarrow+\infty}E_d^{\lambda,\gamma,\delta}\zeta_t(O)=\lim_{t\rightarrow+\infty}E_d^{\lambda,\gamma,\delta}\theta_t(O)=0.
\]
Therefore, by Equations \eqref{equ 3.2} and \eqref{equ 3.3},
\begin{equation}\label{equ 3.5}
\nu^{\lambda,\gamma,\delta}\big(\eta(O)=1\text{~or~}2\big)=0
\end{equation}
when $\lambda<\frac{1}{2d}(1+\frac{1+\delta}{\gamma})$. By Equation \eqref{equ 3.5} and Proposition \ref{proposition 3.2},
\[
\lambda_c(d,\gamma,\delta)\geq \frac{1}{2d}(1+\frac{1+\delta}{\gamma})
\]
and hence $\liminf_{d\rightarrow+\infty}2d\lambda_c(d,\gamma,\delta)\geq 1+\frac{1+\delta}{\gamma}$.

\qed

\section{The proof of $\limsup_{d\rightarrow+\infty}2d\lambda_c(d,\gamma,\delta)\leq 1+\frac{1+\delta}{\gamma}$}\label{section four}
In this section we give the proof of $\limsup_{d\rightarrow+\infty}2d\lambda_c(d,\gamma,\delta)\leq 1+\frac{1+\delta}{\gamma}$. The proof is inspired by the approach introduced in \cite{Xue2017}.

First we introduce a two-stage SIR(susceptible-infected-recovered) model. The two-stage SIR model $\{\rho_t\}_{t\geq 0}$ is a continuous-time Markov process with state space $\{-1,0,1,2\}^{\mathbb{Z}^d}$. The transition rates function of $\{\rho_t\}_{t\geq 0}$ is given as follows. For any $x\in \mathbb{Z}^d$ and $t\geq 0$,
\begin{equation}\label{equ 4.1}
\rho_t(x)\rightarrow i \text{~at rate~}
\begin{cases}
&1 \text{~if~}\rho_t(x)=2\text{~and~}i=-1,\\
&1+\delta \text{~if~}\rho_t(x)=1\text{~and~}i=-1,\\
&\gamma \text{~if~}\rho_t(x)=1\text{~and~}i=2,\\
&\lambda\sum_{y:y\sim x}1_{\{\rho_t(y)=2\}}\text{~if~}\rho_t(x)=0 \text{~and~}i=1,\\
&0 \text{~otherwise},
\end{cases}
\end{equation}
where $\lambda, \gamma, \delta$ is defined as in Equation \eqref{equ 1.1 transition rate}.

Intuitively, for the two-stage SIR model, vertices in state $-1$ are recovered. A recovered vertex can never be infected again and can not infect others. A fully-infected vertex waits for an exponential time with rate $1$ to become recovered while a semi-infected vertex waits for an exponential time with rate $1+\delta$ to become recovered.

Throughout this section we assume that there is no vertex in state $1$ or $-1$ at $t=0$ for the two-stage SIR model. We write $\rho_t$ as $\rho_t^O$ when
$\{x:\rho_0(x)=2\}=\{O\}$. We use $D_t^O$ to denote
\[
\{x:\rho_t^O(x)=2\}
\]
as the set of vertices in state $2$ at the moment $t$ for the two-stage SIR model. We use $P_d^{\lambda, \gamma, \delta}$ to also denote the probability measure
of the two-stage SIR model with parameters $\lambda, \gamma, \delta$. According to the basic coupling of Markov processes (see Section 2.1 of \cite{Lig1985}), it is easy to check that
\begin{equation}\label{equ 4.2}
P_d^{\lambda, \gamma, \delta}\big(C_t^O\neq \emptyset \text{~for all~}t\geq 0\big)\geq P_d^{\lambda, \gamma, \delta}\big(D_t^O\neq \emptyset \text{~for all~}t\geq 0\big)
\end{equation}
for any $\lambda, \gamma, \delta>0$.

For later use, we introduce some independent exponential times. For each $x\in \mathbb{Z}^d$, let $W(x)$ be an exponential time with rate $1$, $Y(x)$ be an exponential time with rate $1+\delta$ while $\Gamma(x)$ be an exponential time with rate $\gamma$. For each pair of neighbors $x,y\in \mathbb{Z}^d$, let $U(x,y)$ be an exponential time with rate $\lambda$. Note that we care about the order of $x$ and $y$, hence $U(x,y)\neq U(y,x)$. We assume that all these exponential times are independent.

For each $n\geq 1$, we define
\begin{align*}
L_n=\Big\{\vec{x}=(x_0,x_1,\ldots,x_n)\in \{\mathbb{Z}^d\}^{n+1}: &x_0=O, x_{i+1}\sim x_i\text{~for all~}0\leq i\leq n-1 \\
&\text{~and~}x_i\neq x_j \text{~for any~}0\leq i<j\leq n\Big\}
\end{align*}
as the set of self-avoiding paths starting at $O$ with length $n$. For each $n\geq 1$ and each $\vec{x}=(x_0,\ldots,x_n)\in L_n$, we use $A_{\vec{x}}$ to denote the event that
\[
U(x_i,x_{i+1})<W(x_i) \text{~and~} \Gamma(x_{i+1})<Y(x_{i+1})
\]
for all $0\leq i\leq n-1$. Then, the two-stage SIR model and these exponential times have the following coupling relationship.
\begin{lemma}\label{lemma 4.1}
$\{\rho_t^O\}_{t\geq 0}$ and $\{W(x)\}_{x\in \mathbb{Z}^d}, \{Y(x)\}_{x\in \mathbb{Z}^d}, \{\Gamma(x)\}_{x\in \mathbb{Z}^d}, \{U(x,y)\}_{x\sim y}$ can be coupled under a same probability space such that for each $n\geq 1$ and any $\vec{x}=(x_0,\ldots,x_n)\in L_n$,
\[
A_{\vec{x}}\subseteq \{x_n\in D_t^O \text{~for some~}t\geq 0\}.
\]
\end{lemma}
According to Lemma \ref{lemma 4.1}, in the sense of coupling, the ending vertex $x_n$ of the self-avoiding path $\vec{x}$ has ever been fully-infected on the event $A_{\vec{x}}$. The detailed proof of Lemma \ref{lemma 4.1} is a little tedious. Here we give an intuitive explanation which is enough to convince Lemma \ref{lemma 4.1}.

\proof[Explanation of Lemma \ref{lemma 4.1}]

The meanings of the exponential times we introduce are as follows. If a vertex $x$ becomes semi-infected at some moment, then $x$ waits for $Y(x)$ units of time to become recovered or waits for $\Gamma(x)$ units of time to become fully-infected, depending on whether $Y(x)<\Gamma(x)$ or $\Gamma(x)<Y(x)$. If $x$ becomes fully-infected at some moment, then $x$ waits for $W(x)$ units of time to become recovered. For any $y\sim x$, the fully-infected vertex $x$ waits for $U(x,y)$ units of time to infect $y$. This infection, which makes $y$ semi-infected, really occurs when and only when $y$ has not been infected by others at an earlier moment and $U(x,y)<W(x)$.

On the event $A_{\vec{x}}$, we can deduce that $x_1,\ldots,x_n$ all belong to $\bigcup_{t\geq 0}D_t^O$ according to the following analysis. For $x_1$, there are two cases. The first case is that $x_1$ is in state $0$ at the moment before $t=U(O,x_1)$, then $x_1$ becomes semi-infected at $t=U(O,x_1)$ since $U(O,x_1)<W(O)$ and $\rho_0(O)=2$. Then, $x_1$ becomes fully-infected at the moment $U(O,x_1)+\Gamma(x_1)$ since $\Gamma(x_1)<Y(x_1)$. The second case is that $x_1$ becomes semi-infected at some moment $s<U(O,x_1)$, then $x$ becomes fully-infected at $s+\Gamma(x_1)$ since $\Gamma(x_1)<Y(x_1)$. In both cases, $x_1$ has ever been fully-infected, i.e., $x_1\in \bigcup_{t\geq 0}D_t^O$. Repeated utilizing of this analysis shows that $x_2, x_3,\ldots,x_n\in \bigcup_{t\geq 0}D_t^O$.

\qed

Inspired by \cite{Xue2017}, we consider a special type of self-avoiding paths. For each $n\geq 1$, we define
\begin{align*}
R_n=\Big\{\vec{x}&=(x_0,\ldots,x_n)\in L_n:x_{i+1}-x_i\in \{\pm e_j:1\leq j\leq d-\lfloor\frac{d}{\log d}\rfloor\}\text{~for any~}i\\
&\text{~such that~}\lfloor\log d\rfloor \nmid(i+1); x_{i+1}-x_i\in \{e_j:d-\lfloor \frac{d}{\log d}\rfloor+1\leq j\leq d\}\\
&\text{~for any~}i\text{~such that~}\lfloor\log d\rfloor \mid(i+1)\Big\},
\end{align*}
where we use $a\mid b$ to denote that $b$ is divisible by $a$ and $\{e_j\}_{1\leq j\leq d}$ are the elementary unit vectors on $\mathbb{Z}^d$, i.e.,
\[
e_j=(0,\ldots,0,\mathop 1\limits_{j \text{th}},0,\ldots,0)
\]
for $1\leq j\leq d$.

According to Lemma \ref{lemma 4.1}, on the event $\bigcap_{n\geq 1}\bigcup_{\vec{x}\in R_n}A_{\vec{x}}$, there are vertices with arbitrarily large $l_1$-norm that have ever been fully-infected, which makes fully-infected vertices survival since each fully-infected vertex waits for an exponential time with rate $1$ to become recovered. Then, by Equation \eqref{equ 4.2},
\begin{align}\label{equ 4.3}
&P^{\lambda,\gamma,\delta}_d\big(C_t^O\neq \emptyset\text{~for all~}t\geq 0\big)\geq P^{\lambda, \gamma, \delta}\big(D_t^O\neq \emptyset \text{~for all~}t\geq 0\big)\notag\\
&\geq P\big(\bigcap_{n\geq 1}\bigcup_{\vec{x}\in R_n}A_{\vec{x}}\big)\geq \lim_{n\rightarrow+\infty}P\big(\bigcup_{\vec{x}\in R_n}A_{\vec{x}}\big).
\end{align}
To bound $P\big(\bigcup_{\vec{x}\in R_n}A_{\vec{x}}\big)$ from below, we introduce a self-avoiding random walk $\{S_n\}_{n\geq 0}$ on $\mathbb{Z}^d$ such that
\[
(S_0,S_1,\ldots,S_n)\in R_n
\]
for each $n\geq 1$. Note that from now on we assume that $d$ is sufficiently large that
\[
2(d-\lfloor \frac{d}{\log d}\rfloor)-\lfloor \log d\rfloor\geq 1.
\]
We define $S_0=O$. For $i\geq 1$ that $\lfloor \log d\rfloor\mid i$,
\[
P\Big(S_i=S_{i-1}+e_l\Big|S_j, 0\leq j\leq i-1\Big)=\frac{1}{\lfloor\frac{d}{\log d}\rfloor}
\]
for each $d-\lfloor\frac{d}{\log d}\rfloor+1\leq l\leq d$. For $i\geq 1$ that $\lfloor \log d\rfloor\nmid i$,
\[
P\Big(S_i=y\Big|S_j,0\leq j\leq i-1\Big)=\frac{1}{|H_{i-1}|}
\]
for any $y\in H_{i-1}$, where
\[
H_{i-1}=\Big\{z:z-S_{i-1}\in \{\pm e_j:1\leq j\leq d-\lfloor\frac{d}{\log d}\rfloor\}\text{~and~}S_j\neq z \text{~for all~}0\leq j\leq i-1\Big\}
\]
 while $|A|$ is the cardinality of the set $A$. Note that $H_{i-1}$ is a random set measurable with respect to the $\sigma$-field generated by $S_0,S_1,\ldots,S_{i-1}$. We claim that
 \begin{equation}\label{equ 4.4}
 |H_{i-1}|\geq 2(d-\lfloor\frac{d}{\log d}\rfloor)-\lfloor\log d\rfloor
 \end{equation}
for each $i\geq 1$. This claim holds according to the following analysis. For each $x=(x_1,\ldots,x_d)\in \mathbb{Z}^d$, we define
\[
u(x)=\sum_{i=d-\lfloor\frac{d}{\log d}\rfloor+1}^d |x_i|,
\]
then $u(S_\cdot)$ increases by $1$ every $\lfloor \log d\rfloor$ steps and hence
\[
\Big|\big\{0\leq j\leq i-1:u(S_j)=u(S_{i-1})\big\}\Big|\leq \lfloor \log d\rfloor.
\]
As a result,
\begin{equation}\label{equ 4.5}
|\{S_0,S_1,\ldots,S_{i-1}\}\bigcap \{S_{i-1}\pm e_j:1\leq j\leq d-\lfloor\frac{d}{\log d}\rfloor\}|\leq \lfloor \log d\rfloor,
\end{equation}
since $u(z)=u(S_{i-1})$ for any $z\in \{S_{i-1}\pm e_j:1\leq j\leq d-\lfloor\frac{d}{\log d}\rfloor\}$. Equation \eqref{equ 4.4} follows from Equation \eqref{equ 4.5} directly.

According to the definition of $\{S_n\}_{n\geq 0}$, it is easy to check that $(S_0,S_1,\ldots,S_n)\in R_n$ for each $n\geq 1$. We let $\{V_n\}_{n\geq 0}$ be an independent copy of $\{S_n\}_{n\geq 0}$ with $V_0=O$. For simplicity, we use $\vec{S}_n$ to denote $(S_0,\ldots,S_n)$ and use $\vec{V}_n$ to denote $(V_0,\ldots, V_n)$ for each $n\geq 1$, then $\vec{S}_n, \vec{V}_n\in R_n$. For any $\vec{x}=(x_0,\ldots,x_n),\vec{y}=(y_0,\ldots,y_n)\in R_n$, we define
\[
F(\vec{x},\vec{y})=\big\{0\leq i\leq n: y_i=x_j \text{~for some~}0\leq j\leq n\big\}
\]
and
\[
K(\vec{x},\vec{y})=\big\{0\leq i\leq n-1: y_i=x_j \text{~and~}y_{i+1}=x_{j+1} \text{~for some~}0\leq j\leq n-1\big\}.
\]
We use $\widehat{P}$ to denote the probability measure of $\{S_n, V_n\}_{n\geq 0}$ and use $\widehat{E}$ to denote the expectation operator with respect to $\widehat{P}$, then the following lemma is crucial for us to prove $\limsup_{d\rightarrow+\infty}2d\lambda_c(d,\gamma,\delta)\leq 1+\frac{1+\delta}{\gamma}$.
\begin{lemma}\label{lemma 4.2}
For each $n\geq 1$,
\[
P\big(\bigcup_{\vec{x}\in R_n}A_{\vec{x}}\big)\geq \frac{1} {\widehat{E}\Big[2^{^{|F(\vec{S}_n,\vec{V}_n)\setminus K(\vec{S}_n,\vec{V}_n)|}}\big(\frac{1+\gamma+\delta}{\gamma}\big)^{|F(\vec{S}_n,\vec{V}_n)|-1}\big(\frac{\lambda+1}{\lambda}\big)^{|K(\vec{S}_n,\vec{V}_n)|}\Big]}.
\]
\end{lemma}

We give the proof of Lemma \ref{lemma 4.2} at the end of this section. Now we show that how to utilize Lemma \ref{lemma 4.2} to prove $\limsup_{d\rightarrow+\infty}2d\lambda_c(d,\gamma,\delta)\leq 1+\frac{1+\delta}{\gamma}$.

\proof[Proof of $\limsup_{d\rightarrow+\infty}2d\lambda_c(d,\gamma,\delta)\leq 1+\frac{1+\delta}{\gamma}$]

Let
\[
F(\vec{S},\vec{V})=\big\{i\geq 0: V_i=S_j \text{~for some~}j\geq 0\big\}
\]
and
\[
K(\vec{S},\vec{V})=\big\{i\geq 0: V_i=S_j \text{~and~}V_{i+1}=S_{j+1} \text{~for some~}j\geq 0\big\},
\]
then
\[
\lim_{n\rightarrow+\infty}|K(\vec{S}_n,\vec{V}_n)|=|K(\vec{S},\vec{V})| \text{~while~}\lim_{n\rightarrow+\infty}|F(\vec{S}_n,\vec{V}_n)|=|F(\vec{S},\vec{V})|
\]
and hence
\begin{align*}
&\lim_{n\rightarrow+\infty}\widehat{E}\Big[2^{^{|F(\vec{S}_n,\vec{V}_n)\setminus K(\vec{S}_n,\vec{V}_n)|}}\big(\frac{1+\gamma+\delta}{\gamma}\big)^{|F(\vec{S}_n,\vec{V}_n)|-1}\big(\frac{\lambda+1}{\lambda}\big)^{|K(\vec{S}_n,\vec{V}_n)|}\Big]  \\
&=\widehat{E}\Big[2^{^{|F(\vec{S},\vec{V})\setminus K(\vec{S},\vec{V})|}}\big(\frac{1+\gamma+\delta}{\gamma}\big)^{|F(\vec{S},\vec{V})|-1}\big(\frac{\lambda+1}{\lambda}\big)^{|K(\vec{S},\vec{V})|}\Big]\\
&=\big(\frac{\gamma}{1+\gamma+\delta}\big)\widehat{E}\Big[\big(\frac{2(1+\gamma+\delta)}{\gamma}\big)^{|F(\vec{S},\vec{V})\setminus K(\vec{S},\vec{V})|}\big(\frac{\lambda+1}{\lambda}\frac{1+\gamma+\delta}{\gamma}\big)^{|K(\vec{S},\vec{V})|}\Big]
\end{align*}
according to Dominated Convergence Theorem. Then by Equation \eqref{equ 4.3} and Lemma \ref{lemma 4.2},
\begin{align}\label{equ 4.6}
&P^{\lambda,\gamma,\delta}_d\big(C_t^O\neq \emptyset\text{~for all~}t\geq 0\big) \\
&\geq \frac{1}{\big(\frac{\gamma}{1+\gamma+\delta}\big)\widehat{E}\Big[\big(\frac{2(1+\gamma+\delta)}{\gamma}\big)^{|F(\vec{S},\vec{V})\setminus K(\vec{S},\vec{V})|}\big(\frac{\lambda+1}{\lambda}\frac{1+\gamma+\delta}{\gamma}\big)^{|K(\vec{S},\vec{V})|}\Big]}. \notag
\end{align}
Reference \cite{Xue2017} gives a detailed calculation of the function
\[
J(C_1,C_2)=\widehat{E}\Big[C_1^{|F(\vec{S},\vec{V})\setminus K(\vec{S},\vec{V})|}C_2^{|K(\vec{S},\vec{V})|}\Big].
\]
According to Lemma 3.4 of \cite{Xue2017}, there exists $M_1, M_2$ which do not depend on $C_1, C_2$ and the dimension $d$ of the lattice that
\[
J(C_1, C_2)\leq M_2C_1\sum_{n=0}^{+\infty}\Big[\frac{(\lfloor\log d\rfloor^{^{\frac{3}{\lfloor\log d\rfloor-1}}})C_2}{2(d-\lfloor\frac{d}{\log d}\rfloor)-\lfloor\log d\rfloor}+\frac{C_2}{\lfloor\frac{d}{\log d}\rfloor\lfloor\log d\rfloor^3}+\frac{M_1(\log d)^5C_1}{d}\Big]^n
\]
for any $C_1, C_2>0$.

For given $\vartheta>1$, let $\lambda=\frac{\vartheta}{2d}\frac{1+\gamma+\delta}{\gamma}$, then it is easy to check that
\[
\frac{(\lfloor\log d\rfloor^{^{\frac{3}{\lfloor\log d\rfloor-1}}})\frac{\lambda+1}{\lambda}\frac{1+\gamma+\delta}{\gamma}}{2(d-\lfloor\frac{d}{\log d}\rfloor)-\lfloor\log d\rfloor}+\frac{\frac{\lambda+1}{\lambda}\frac{1+\gamma+\delta}{\gamma}}{\lfloor\frac{d}{\log d}\rfloor\lfloor\log d\rfloor^3}+\frac{M_1(\log d)^5\frac{2(1+\gamma+\delta)}{\gamma}}{d}<1
\]
for sufficiently large $d$ and hence
\[
J\big(\frac{2(1+\gamma+\delta)}{\gamma}, \frac{\lambda+1}{\lambda}\frac{1+\gamma+\delta}{\gamma}\big)<+\infty
\]
for sufficiently large $d$. As a result, by Equation \eqref{equ 4.6},
\[
P^{\lambda,\gamma,\delta}_d\big(C_t^O\neq \emptyset\text{~for all~}t\geq 0\big)\geq \frac{1}{\big(\frac{\gamma}{1+\gamma+\delta}\big)J\big(\frac{2(1+\gamma+\delta)}{\gamma}, \frac{\lambda+1}{\lambda}\frac{1+\gamma+\delta}{\gamma}\big)}>0
\]
when $\lambda=\frac{\vartheta}{2d}\frac{1+\gamma+\delta}{\gamma}$ and $d$ is sufficiently large. Therefore,
\[
\lambda_c(d,\gamma,\delta)\leq \frac{\vartheta}{2d}\frac{1+\gamma+\delta}{\gamma}=\frac{\vartheta}{2d}(1+\frac{1+\delta}{\gamma})
\]
for sufficiently large $d$ and hence
\[
\limsup_{d\rightarrow+\infty}2d\lambda_c(d,\gamma,\delta)\leq \vartheta(1+\frac{1+\delta}{\gamma}).
\]
Since $\vartheta>1$ is arbitrary, let $\vartheta \downarrow 1$ then the proof is complete.

\qed

To finish this section, we need to prove Lemma \ref{lemma 4.2}. The proof of Lemma \ref{lemma 4.2} relies heavily on the following proposition, which is Lemma 3.3 of \cite{Xue2017}.
\begin{proposition}\label{proposition 4.3}(Xue, 2017)
If $B_1, B_2,\ldots, B_n$ are $n$ arbitrary events defined under the same probability space such that $P(B_i)>0$ for $1\leq i\leq n$ and $p_1,p_2,\ldots,p_n$ are $n$ positive constants such that
$\sum_{j=1}^np_j=1$, then
\[
P(\bigcup_{j=1}^{+\infty}B_j)\geq \frac{1}{\sum\limits_{i=1}^n\sum\limits_{j=1}^np_ip_j\frac{P(B_i\bigcap B_j)}{P(B_i)P(B_j)}}.
\]
\end{proposition}
For the proof of Proposition \ref{proposition 4.3}, see Section 3 of \cite{Xue2017}. At last we give the proof of Lemma \ref{lemma 4.2}.

\proof[Proof of Lemma \ref{lemma 4.2}]

For each $\vec{x}\in R_n$, let $p_{\vec{x}}$ be the probability that $\vec{S}_n=\vec{x}$, then by Proposition \ref{proposition 4.3},
\begin{equation}\label{equ 4.7}
P\big(\bigcup_{\vec{x}\in R_n}A_{\vec{x}}\big)\geq \frac{1}{\sum_{\vec{x}\in R_n}\sum_{\vec{y}\in R_n}p_{\vec{x}}p_{\vec{y}}\frac{P\big(A_{\vec{x}}\bigcap A_{\vec{y}}\big)}{P\big(A_{\vec{x}}\big)P\big(A_{\vec{y}}\big)}}.
\end{equation}
Now we bound $\frac{P\big(A_{\vec{x}}\bigcap A_{\vec{y}}\big)}{P\big(A_{\vec{x}}\big)P\big(A_{\vec{y}}\big)}$ from above. For
\[
\vec{x}=(x_0,\ldots, x_n), \vec{y}=(y_0, \ldots, y_n)\in R_n,
\]
if $x_i\not\in \big\{y_0,\ldots,y_n\big\}$ for some $0<i\leq n$, then the factor
\[
P\big(U(x_i,x_{i+1})<W(x_i)\big)P\big(\Gamma(x_i)<Y(x_i)\big)
\]
appears once in both $P\big(A_{\vec{x}}\bigcap A_{\vec{y}}\big)$ and $P\big(A_{\vec{x}}\big)P\big(A_{\vec{y}}\big)$, which can be cancelled. Similarly, if $y_j\not\in \big\{x_0,\ldots,x_n\big\}$ for some $0<j\leq n$, then the factor
\[
P\big(U(y_j,y_{j+1})<W(y_j)\big)P\big(\Gamma(y_j)<Y(y_j)\big)
\]
appears once in both $P\big(A_{\vec{x}}\bigcap A_{\vec{y}}\big)$ and $P\big(A_{\vec{x}}\big)P\big(A_{\vec{y}}\big)$, which can be cancelled. If $j\in F(\vec{x},\vec{y})\setminus \{0\}$, then $x_i=y_j$ for some $0<i\leq n$ and the factor
\[
P\big(\Gamma(y_j)<Y(y_j)\big)=P\big(\Gamma(x_i)<Y(x_i)\big)=\frac{\gamma}{1+\gamma+\delta}
\]
appears twice in $P\big(A_{\vec{x}}\big)P\big(A_{\vec{y}}\big)$ but appears once in $P\big(A_{\vec{x}}\bigcap A_{\vec{y}}\big)$ since
\[
\{\Gamma(y_j)<Y(y_j)\}\bigcap \{\Gamma(x_i)<Y(x_i)\}=\{\Gamma(y_j)<Y(y_j)\},
\]
which generates a factor $\frac{1+\gamma+\delta}{\gamma}$ for $\frac{P\big(A_{\vec{x}}\bigcap A_{\vec{y}}\big)}{P\big(A_{\vec{x}}\big)P\big(A_{\vec{y}}\big)}$.
If $j\in K(\vec{x},\vec{y})$, then $x_i=y_j$ while $x_{i+1}=y_{j+1}$ for some $0\leq i\leq n-1$ and the factor
\[
P\big(U(y_j, y_{j+1})<W(y_j)\big)=P\big(U(x_i,x_{i+1})<W(x_i)\big)=\frac{\lambda}{1+\lambda}
\]
appears twice in $P\big(A_{\vec{x}}\big)P\big(A_{\vec{y}}\big)$ but appears once in $P\big(A_{\vec{x}}\bigcap A_{\vec{y}}\big)$ since
\[
\{U(y_j, y_{j+1})<W(y_j)\}\bigcap\{U(x_i,x_{i+1})<W(x_i)\}=\{U(y_j, y_{j+1})<W(y_j)\},
\]
which generates a factor $\frac{\lambda+1}{\lambda}$ for $\frac{P\big(A_{\vec{x}}\bigcap A_{\vec{y}}\big)}{P\big(A_{\vec{x}}\big)P\big(A_{\vec{y}}\big)}$. If $j\in F(\vec{x},\vec{y})\setminus K(\vec{x},\vec{y})$, then $x_i=y_j$ while $x_{i+1}\neq y_{j+1}$ for some $0\leq i\leq n-1$ and hence $P\big(A_{\vec{x}}\bigcap A_{\vec{y}}\big)$ has the factor $P\big(U(x_i,x_{i+1})<W(x_i), U(x_i,y_{j+1})<W(x_i)\big)$ while $P\big(A_{\vec{x}}\big)P\big(A_{\vec{y}}\big)$ has the factor $P\big(U(x_i,x_{i+1})<W(x_i)\big)P\big(U(x_i,y_{j+1})<W(x_i)\big)$, which generates a factor
\[
\frac{P\big(U(x_i,x_{i+1})<W(x_i), U(x_i,y_{j+1})<W(x_i)\big)}{P\big(U(x_i,x_{i+1})<W(x_i)\big)P\big(U(x_i,y_{j+1})<W(x_i)\big)}
=\frac{2\lambda+2}{2\lambda+1}\leq 2
\]
for $\frac{P\big(A_{\vec{x}}\bigcap A_{\vec{y}}\big)}{P\big(A_{\vec{x}}\big)P\big(A_{\vec{y}}\big)}$. In conclusion,
\begin{align}\label{equ 4.8}
\frac{P\big(A_{\vec{x}}\bigcap A_{\vec{y}}\big)}{P\big(A_{\vec{x}}\big)P\big(A_{\vec{y}}\big)}&\leq \big(\frac{1+\gamma+\delta}{\gamma}\big)^{|F(\vec{x}_n,\vec{y}_n)\setminus\{O\}|}\big(\frac{1+\lambda}{\lambda}\big)^{|K(\vec{x},\vec{y})|}\big(2\big)^{^{|F(\vec{x},\vec{y})\setminus K(\vec{x},\vec{y})|}}  \notag\\
&=2^{^{|F(\vec{x},\vec{y})\setminus K(\vec{x},\vec{y})|}}\big(\frac{1+\gamma+\delta}{\gamma}\big)^{|F(\vec{x},\vec{y})|-1}\big(\frac{\lambda+1}{\lambda}\big)^{|K(\vec{x},\vec{y})|}.
\end{align}
By Equation \eqref{equ 4.8} and the definition of $p_{\vec{x}}$,
\begin{align}\label{equ 4.9}
&\sum_{\vec{x}\in R_n}\sum_{\vec{y}\in R_n}p_{\vec{x}}p_{\vec{y}}\frac{P\big(A_{\vec{x}}\bigcap A_{\vec{y}}\big)}{P\big(A_{\vec{x}}\big)P\big(A_{\vec{y}}\big)}\notag\\
&\leq \sum_{\vec{x}\in R_n}\sum_{\vec{y}\in R_n}p_{\vec{x}}p_{\vec{y}}\Big[2^{^{|F(\vec{x},\vec{y})\setminus K(\vec{x},\vec{y})|}}\big(\frac{1+\gamma+\delta}{\gamma}\big)^{|F(\vec{x},\vec{y})|-1}\big(\frac{\lambda+1}{\lambda}\big)^{|K(\vec{x},\vec{y})|}\Big] \notag\\
&=\widehat{E}\Big[2^{^{|F(\vec{S}_n,\vec{V}_n)\setminus K(\vec{S}_n,\vec{V}_n)|}}\big(\frac{1+\gamma+\delta}{\gamma}\big)^{|F(\vec{S}_n,\vec{V}_n)|-1}\big(\frac{\lambda+1}{\lambda}\big)^{|K(\vec{S}_n,\vec{V}_n)|}\Big].
\end{align}
Lemma \ref{lemma 4.2} follows from Equations \eqref{equ 4.7} and \eqref{equ 4.9} directly. 

\qed

\quad

\textbf{Acknowledgments.} The author is grateful to the financial
support from the National Natural Science Foundation of China with
grant number 11501542 and the financial support from Beijing
Jiaotong University with grant number KSRC16006536.

{}
\end{document}